\documentclass[final,5p,twocolumn]{elsarticle}

\usepackage{graphicx}          % Include this line if your
\usepackage[dvips]{epsfig}    % or this line, depending on which
\usepackage{xcolor}

\usepackage{amssymb}
\usepackage{amsthm}

\usepackage{amsmath}          % Include this line if your

\usepackage{algorithmic}

\journal{Communications in Nonlinear Science and Numerical Simulation}

\begin{document}

\begin{frontmatter}

\title{On convergence analysis of feedback control with integral
action and discontinuous relay perturbation\tnoteref{label1}} % Title, preferably not more

\author{Michael Ruderman}
\ead{michael.ruderman@uia.no}

%\fntext[*]{AAA aaa bbb }

\tnotetext[label1]{\textcolor[rgb]{0.00,0.07,1.00}{AUTHOR'S ACCEPTED MANUSCRIPT}}


\address{University of Agder, Department of Engineering Sciences, Norway}

\begin{abstract}                          % Abstract of not more than 200 words.
We consider third-order dynamic systems which have an integral
feedback action and discontinuous relay disturbance. More
specifically for the applications, the focus is on the integral
plus state-feedback control of the motion systems with
discontinuous Coulomb-type friction. We recall the stiction region
is globally attractive where the resulting hybrid system has also
solutions in Filippov sense, while the motion trajectories remain
in that idle state (called in tribology as stiction) until the
formulated sliding-mode condition is violated by the growing
integral feedback quantity. We analyze the conditions for
occurrence of the slowly converging stick-slip cycles. We also
show that the hybrid system is globally but only asymptotically
stable, and almost always not exponentially. A particular case of
the exponential convergence can appear for some initial values,
assuming the characteristic equation of the linear subsystem has
dominant real roots. Illustrative numerical examples are provided
alongside with the developed analysis. In addition, a laboratory example
is shown with experimental evidence to support the convergence
analysis provided.
\end{abstract}

\begin{keyword}
Nonlinear system \sep relay feedback \sep Coulomb friction \sep
convergence analysis  \sep hybrid dynamics
\end{keyword}

\end{frontmatter}

\newtheorem{thm}{Theorem}
\newtheorem{lem}[thm]{Lemma}
\newtheorem{clr}{Corollary}
\newdefinition{rmk}{Remark}
\newproof{pf}{Proof}

%%%%%%%%%%%%%%%%%%%%%%%%%%%%%%%%%%%%%%%%%%%%%%%%%%%%%%%%%%%%%%%%%%%%%%%%%%%%%%%%
\section{Introduction}
\label{sec:1}

Discontinuous relays in feedback are frequently appearing in the
dynamic systems leading, thus, to a hybrid or switching behavior
see e.g. \cite{liberzon2003} and, as particular case, to the
sliding modes, see e.g. \cite{shtessel2014}. The latter usually
require the solutions in Filippov sense (see \cite{filippov1988}
for detail) of the differential equations with discontinuous
right-hand-side -- the methodology of which we will also make use
in the present study. Several works considered extensively both,
the existence of fast switches and the limit cycles in relay
feedback systems, see \cite{johansson1999,johansson2002}. It has
been discussed in \cite{johansson2002}, that `if the linear part
of the relay feedback system has pole excess one and certain other
conditions are fulfilled, then the system has a limit cycle with
sliding mode.' At the same time, a very common application example
of the integral plus state-feedback control in the motion systems
with discontinuous Coulomb-type friction has the pole excess one,
while being usually understood as globally asymptotically stable
for the properly assigned linear feedback gains, see e.g.
\cite{armstrong1996,bisoffi2017,ruderman2022}.

The complex and rather parasitic (for applications) effects of the
Coulomb friction in context of a controlled motion are since long
studied in system and control communities, see e.g.
\cite{armstrong1994}. The summarized basics of the kinetic
friction and how it affects the feedback controlled systems can be
found, for instance, in the recent text \cite{ruderman2023}. Also
the experimental studies, which are disclosing the impact of
Coulomb friction in vicinity to the reference settling point, e.g.
\cite{ruderman2016}, can be consulted in addition. How a feedback
controlled second-order system (i.e. without integral control
action) comes to the so-called stiction due to the relay-type
Coulomb friction was analyzed in \cite{alvarez2000}. A possible
appearance of the friction-induced limit cycles was addressed in
e.g. \cite{radcliffe1990,olsson2001}. The source of such
persistent limit cycles is, however, not trivial and requires more
than only a relay-type Coulomb nonlinearity for mapping the
frictional interactions. However, a simple relay-type Coulomb
friction force provides already a sufficient basis for studying
the degradation of feedback control and the corresponding
shortcomings of an integral control action. According to the
focused studies \cite{armstrong1996,bisoffi2017,ruderman2022}, the
global asymptotic stability (GAS) of a
proportional-integral-derivative (PID) controlled system with
relay-type Coulomb friction is given, and the slowly converging
stick-slip cycles appear instead of persistent limit cycles.

Against the background mentioned above, the present work aims to
make a further contribution to the convergence analysis of
feedback control with an integral action, when a discontinuous
relay occurs in feedback as a known perturbation. Below, the
problem statement is formulated, followed by the main results
placed in section \ref{sec:2}. The
dedicated illustrative numerical examples are provided in section
\ref{sec:3}. In order to reinforce the main statements and
conclusions of the analysis performed, a laboratory experimental
example is also provided qualitatively, in section \ref{sec:4}.
Section \ref{sec:5} summarizes the paper and highlights several
points for discussion. The particular solutions in use are taken
into Appendix.

\subsection*{Problem statement}
\label{sec:1:sub:1}

We consider the following class of third-order systems
\begin{equation}
\dddot{y}(t) + a \, \ddot{y}(t) + b \, \dot{y}(t) + c \, y(t) = -
\gamma \, \mathrm{sign}\bigl(\ddot{y}(t)\bigr), \label{eq:ode}
\end{equation}
where $a,\, b,\, c > 0$ are the design parameters, and $\gamma >
0$ is the gain of the discontinuous relay function assumed to be
known. The latter, described by the sign operator, is defined in
the Filippov sense, see \cite{filippov1988}, on the closed
interval in zero, i.e.
\begin{equation}
\label{eq:relay}
    \mathrm{sign}(z)= \left\{%
\begin{array}{ll}
    1, & \; z>0, \\
    \left[-1, 1\right], & \; z=0,\\
    -1, & \; z<0. \\
\end{array}%
\right.
\end{equation}

%\vspace{2mm}

\begin{rmk}
\label{rem:1} The system class \eqref{eq:ode} can well describe
the second-order state-feedback controlled motion with
discontinuous Coulomb friction, cf.
\cite{ruderman2022,ruderman2023}, and an additional integral
output feedback. This case, the constants $a$ and $b$ capture
simultaneously both, the state feedback gains and the mechanical
parameters of the system. Same time, $c$ and $\gamma$ constitute
the integral feedback gain and the Coulomb friction coefficient,
respectively. Moreover, all parameters are normalized by a
positive inertial mass.
\end{rmk}

%\vspace{2mm}

Further we note that generally
\begin{equation}
y(0) = C_1, \quad \dot{y}(0) = C_2, \quad \ddot{y}(0) = C_3,
\label{eq:initcond}
\end{equation}
where $C_1, C_2, C_3 \in \mathbb{R}$, while for a typical motion
control the scenario with $C_1 = C_3 = 0$ and $C_2 \neq 0$ is very
common. The objectives targeted in the following are to
\begin{enumerate}[(i)]
    \item provide a esuriently simple proof of the global asymptotic stability (GAS),
cf. with \cite{bisoffi2017,ruderman2022}, of the system
\eqref{eq:ode} and show the parametric conditions that guarantee
the GAS property;
    \item given (i), to show where in the state-space the system
    \eqref{eq:ode} becomes transiently sticking, i.e.
    $\ddot{y}(t) = 0$ for all $t \in [t_{i}^{s}, \,t_{i}^{c}]$ with
    $0 \leq t_{i}^{s} < t_{i}^{c} < \infty$ and $i \in
    \mathbb{N}^{+}$, and how it depends on the system parameters;
    \item given (i) and (ii), to analyze the asymptotic convergence
    in the periodic stick-slip cycles, while a slipping phase is
    characterized by $\ddot{y}(t) \neq 0$ for all $t \in [t_{i}^{c},
    \,t_{i+1}^{s}]$ with $0 \leq t_{i}^{c} < t_{i+1}^{s} <
    \infty$ and $i \in \mathbb{N}^{+}$.
\end{enumerate}

%%%%%%%%%%%%%%%%%%%%%%%%%%%%%%%%%%%%%%%%%%%%%%%%%%%%%%%%%%%%%%%%%%%%%%%%%%%%%%%%
\section{Main results} \label{sec:2}

\subsection{Global asymptotic stability} \label{sec:2:sub:1}

\begin{figure}[!h]
\centering
\includegraphics[width=0.8\columnwidth]{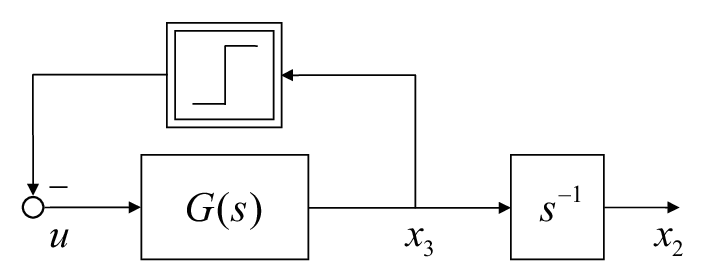}
\caption{Structure of the third-order system with feedback relay
perturbation.} \label{fig:1}
\end{figure}

From the structural viewpoint, the system \eqref{eq:ode} can be
represented as feedback loop with the discontinuous relay
\eqref{eq:relay} as shown in Figure \ref{fig:1}. The linear
subsystem described by the
transfer function
\begin{equation}
G(s) = \frac{y(s) s^2}{u(s)} = C (sI - A)^{-1} B \label{eq:tf}
\end{equation}
is of third-order with
\begin{equation}\label{eq:ssmatr}
   A= \left(%
\begin{array}{ccc}
  0 & 1 & 0 \\
  0 & 0 & 1 \\
  -c & -b & -a \\
\end{array}%
\right),
B=\left(%
\begin{array}{c}
  0 \\
  0 \\
  \gamma \\
\end{array}%
\right),
C^{\top}=\left(%
\begin{array}{c}
  0 \\
  0 \\
  1 \\
\end{array}%
\right),
\end{equation}
to be the system matrix and input and output coupling vectors,
respectively. The $3\times3$ identity matrix is denoted by $I$,
and $s$ is the complex Laplace variable. The corresponding state
vector is $x \equiv [x_1,\, x_2,\, x_3]^{\top} = [y,\, \dot{y},\,
\ddot{y}]^{\top}$. Note that the feedback relay acts as an input
perturbation to $G(s)$.
We are now in the position to formulate and prove the GAS
properties of the system \eqref{eq:ode}, \eqref{eq:relay}.

\begin{thm}
\label{theorem:1} The system \eqref{eq:ode}, \eqref{eq:relay} with
$a,\,b,\,c,\,\gamma > 0$ is globally asymptotically stable (GAS)
if
\begin{equation}
ab > c. \label{eq:GAScond}
\end{equation}
\end{thm}

\begin{pf}
Proof is carried out in the three following steps.

\emph{First} -- we demonstrate that between two consecutive
switches of the relay \eqref{eq:relay} the linear subsystem in
\eqref{eq:ode}, i.e. its left-hand-side, is GAS. For the
input-output system \eqref{eq:tf}, cf. Figure \ref{fig:1}, assume
a quadratic Lyapunov function
\begin{equation}
\label{eq:lyapfun} V(x) = x^{\top} P x \quad \hbox{with} \quad
P=P^{\top}.
\end{equation}
Let us use the fact that if $A$ is Hurwitz, i.e.
$$
\mathrm{Re} \{ \lambda_j \} < 0
$$
for all eigenvalues $\lambda_j$ of $A$, then $P$ represents the
unique solution of the Lyapunov equation
\begin{equation}
\label{eq:lyapeq} PA + A^{\top} P = -Q,
\end{equation}
for a given positive definite symmetric matrix $Q$, cf.
\cite[Theorem~4.6]{khalil2002} and
\cite[Lemma~4.28]{antsaklis2007}. In other words, let us find a
positive definite symmetric matrix $P$ satisfying
\eqref{eq:lyapeq}. One can show that assuming $Q = I$, the
Lyapunov equation \eqref{eq:lyapeq} has the analytic solution $P$,
which is given in Appendix \emph{A}. Now, following
\cite[Theorem~4.6]{khalil2002}, let us derive the parametric
conditions for which $A$ is Hurwitz. This will imply that the
determined $P$ is the unique solution of \eqref{eq:lyapeq} and
also positive definite. Applying the algebraic criterion of
Routh-Hurwitz (see e.g. \cite{franklin2020}) to the characteristic
polynomial
\begin{equation}
\label{eq:charpol} \alpha(s) = s^3 + a s^2 + b s + c = 0
\end{equation}
of the system \eqref{eq:tf}, \eqref{eq:ssmatr}, we obtain
\eqref{eq:GAScond}. This proves that the system \eqref{eq:tf},
\eqref{eq:ssmatr} remains always GAS between two consecutive
switchings of the relay at $0 < t_i^0 < t_{i+1}^0 < \infty$ and,
moreover, that the trajectory $x(t)$ continues to converge
exponentially and uniformly on any time interval $[t_i^0, \,
t_{i+1}^0]$.

\emph{Second} -- for the system \eqref{eq:tf}, \eqref{eq:ssmatr},
since it is GAS on an interval $[t_i^0, \, t_{i+1}^0]$, we apply
the final value theorem for the switched input $u(t_i^0) = \pm 2
\gamma \, h(t_i^0)$, where $h(\cdot)$ is Heaviside (i.e. unit)
step function, and obtain for steady-state
\begin{eqnarray}
\label{eq:finalx1}
  x_1(0) &=& \pm \frac{\gamma}{\alpha(s)}
\frac{2s}{s} \, \Bigl|_{s \rightarrow 0} = \pm \frac{2\gamma}{c}, \\
  x_2(0) &=& \pm \frac{\gamma}{\alpha(s)}
\frac{2s^2}{s} \, \Bigl|_{s \rightarrow 0} = 0. \label{eq:finalx2}
\end{eqnarray}
For the corresponding input to the relay, i.e. $x_3$ cf. Figure
\ref{fig:1}, one obtains (in Laplace domain)
\begin{equation}
\label{eq:x3value} x_3(s) = \pm \frac{2 \gamma \,s }{\alpha(s)},
\end{equation}
which is converging asymptotically to zero as $s \rightarrow 0$ if
$t_{i+1}^0 \rightarrow \infty$. This completes the proof of
asymptotic convergence in case there is no zero-crossing of
$x_3(t)$ and, consequently, no subsequent switching at $t_{i+1}^0
< \infty$. Otherwise, if there is a subsequent zero-crossing of
$x_3(t)$:

\emph{Third} -- we should examine whether the switchings occur
permanently in the sequence and, as a result, a non-converging
periodic solution in the sense of a limit cycle occurs. Using the
describing function analysis method, see e.g. \cite{Atherton75}
for details, the harmonic balance equation
$$
G(j\omega) N(\mathcal{A}, \omega) + 1 = 0,
$$
where $N(\cdot)$ is the describing function of a static
nonlinearity in feedback (here the relay) and $\mathcal{A}$ is the
amplitude of the first harmonic of the periodic solution, results
in
\begin{equation}
\label{eq:harmbalance} 1 - \frac{4 \gamma \omega^2}{\pi
\mathcal{A} \bigl( - j \omega^3 - a \omega^2 + j b \omega + c
\bigr)} = 0.
\end{equation}
Solving \eqref{eq:harmbalance} with respect to $\mathcal{A}$ and,
then, evaluating the parameters condition for $\mathcal{A}$ to be
a real number, one can show that \eqref{eq:harmbalance} has no
real solutions in terms of the $\mathcal{A}$ and $\omega$. It
implies that there is no intersection point of the negative
reciprocal of $N(\mathcal{A}) = 4 \gamma \, (\pi
\mathcal{A})^{-1}$ with the Nyquist plot $G(j\omega)$. Therefore,
there is no non-converging periodic solution of the system
\eqref{eq:tf}, \eqref{eq:ssmatr}. This completes the proof.
\end{pf}

\begin{rmk}
\label{rem:2} Note that the GAS condition \eqref{eq:GAScond} is
conservative, since the passive relay in feedback (Figure
\ref{fig:1}) provides an additional system damping, i.e. energy
dissipation during each slipping phase defined by (iii) in section
\ref{sec:1:sub:1}.
\end{rmk}

\begin{rmk}
\label{rem:2:1} Note that the standard absolute stability
criterion, see \cite{khalil2002} for details, cannot be applied
despite the system \eqref{eq:tf}, \eqref{eq:ssmatr} clearly
possesses the structure of a Lur'e problem \cite{lur1944}, which
means the system has a static nonlinearity belonging to the sector
$[0, \infty]$. The non-applicability of the absolute stability
criterion is due to the fact that the linear subsystem
\eqref{eq:tf} is not strictly positive real, cf.
\cite[Theorem~7.1]{khalil2002}.
\end{rmk}

\subsection{Properties of stiction region} \label{sec:2:sub:2}

Despite the system \eqref{eq:ode}, \eqref{eq:relay} proves to be
GAS, provided \eqref{eq:GAScond} is satisfied, the $(x_2,x_3)$
trajectories can temporarily come to stiction, i.e. $x_3 = 0$, due
to the switching relay, cf. \cite{ruderman2022}. Defining the
switching surface as
\begin{equation}
\label{eq:slidingsurf} S : = \bigl\{ x \in \mathbb{R}^3  \: | \: C
x=0 \bigr \},
\end{equation}
the state trajectories can either cross it or they can stay on it
(then in sliding-mode, e.g. \cite{shtessel2014}) depending on the
vector field of the corresponding differential inclusion, cf.
\cite{liberzon2003},
\begin{equation}
\label{eq:inclusion} \dot{x} \in F(x)= Ax - B\, \mathrm{sign}
(x_3),
\end{equation}
in vicinity to $S$. Here, the relay value depends on whether $x
\rightarrow S$ from $x \in \mathbb{R}^3_{+} \equiv \{ \mathbb{R}^3
\backslash S \, : \: x_3 > 0 \}$ or from $x \in \mathbb{R}^3_{-}
\equiv \{ \mathbb{R}^3 \backslash S \, : \: x_3 < 0 \}$. Let us
establish the condition and prove the appearance of the
sliding-mode for \eqref{eq:ode}, \eqref{eq:relay}.

\begin{thm}
\label{theorem:2} Assume the system \eqref{eq:ode},
\eqref{eq:relay} with $a,\,b,\,c,\,\gamma > 0$ is GAS in the sense
of Theorem \ref{theorem:1}. Then it is sticking in terms of
$\ddot{y}=0$ and has an unique solution in Filippov sense if
\begin{equation} \gamma > \bigl| b\dot{y} + cy \bigr|.
\label{eq:theorem2}
\end{equation}
\end{thm}

\begin{pf}
The proof is carried out by using the standard existence condition
of the sliding-modes, cf. \cite{shtessel2014}.

Excluding, for a while, the global attractiveness property (that
will be analyzed later), we use the condition
$$
S \dot{S} < 0 \quad \hbox{for} \quad S \neq 0,
$$
and, more specifically, the \emph{sufficient condition}
\begin{equation}
\underset{S \rightarrow 0^{+}} {\lim} \, \dot{S} < 0 \quad
\hbox{and} \quad \underset{S \rightarrow 0^{-}} {\lim} \, \dot{S}
> 0
\label{eq:suffcond}
\end{equation}
for the sliding-mode to exist, cf. \cite[Chapter~2]{utkin2020}.
Taking the time derivative of $S = Cx$ and evaluating it we obtain
\begin{equation}
\dot{S} = -\bigl(a x_3 + b x_2 + c x_1 + \gamma \,
\mathrm{sign}(x_3)\bigr). \label{eq:Dsurface}
\end{equation}
Substituting \eqref{eq:Dsurface} into \eqref{eq:suffcond} and
evaluating both limits results in
\begin{equation}
\gamma > \left\{%
\begin{array}{ll}
    -(b x_2 + c x_1), & \hbox{ if } S \rightarrow 0^{+}, \\[1mm]
    b x_2 + c x_1, & \hbox{ if } S \rightarrow 0^{-}.
\end{array}%
\right. \label{eq:evalcond}
\end{equation}
The inequality \eqref{eq:evalcond} is equivalent to
\eqref{eq:theorem2} for $S=0$ when $x_3 \equiv \ddot{y} = 0$. This
completes the proof.
\end{pf}

\begin{rmk}
\label{rem:3} Condition \eqref{eq:theorem2} is equivalent to the
set of attraction $\{x\in S:|CAx|<|CB|\}$ if a relay feedback
system\footnote{Note that here $|\cdot|$ denotes the
absolute-value norm.} has $CB>0$. Note that the latter was
demonstrated in \cite[Section~4]{johansson1999} in context of a
necessary and sufficient condition for the systems to have the
consecutive fast relay switchings.
\end{rmk}

\begin{figure}[!h]
\centering
\includegraphics[width=0.8\columnwidth]{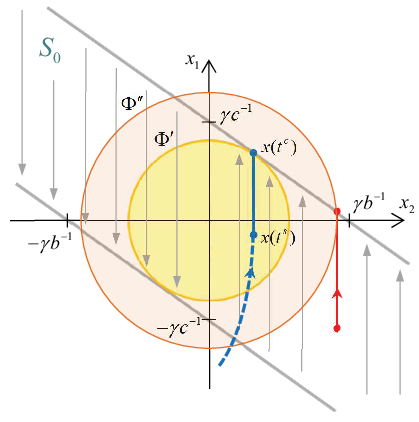}
\caption{Stiction region $S_0$ in $\bigl(x_2, x_1\bigr)$
coordinates.} \label{fig:2}
\end{figure}

The parametric condition \eqref{eq:theorem2} can be well
interpreted graphically as shown in Figure \ref{fig:2}.
The stiction region
\begin{equation}
S_0 \equiv \{x \in \mathbb{R}^3 \, : \: x_3=0, \, | bx_2 + cx_1 |
< \gamma \} \label{eq:sticreg}
\end{equation}
lies in the $(x_2,x_1)$ plane, bounded by two straight lines which
cross the $x_2$ and $x_1$ axes at $\pm \gamma b^{-1}$ and $\pm
\gamma c^{-1}$. During the system is sticking, the vector field
$F$ is:
\begin{enumerate}[(i)]
    \item lying between both lines that satisfies \eqref{eq:sticreg},
    \item orthogonal to the $x_2$-axis,
    \item unambiguously directed towards the growing $x_1$
values for $x_2 > 0$, and towards the falling $x_1$ values for
$x_2 < 0$, see grey arrows shown in Figure \ref{fig:2}.
\end{enumerate}
In the following, the time instants of reaching $S_0$ will be
denoted by $t_i^s$ and the time instants of leaving $S_0$ by
$t_i^c$ for $\forall \; i \in \mathbb{N}^{+}$, cf. Problem
statement in section \ref{sec:1}.

\begin{rmk}
\label{rem:4} Once reaching $S_0$ at some $0 \leq t^s < \infty$,
the system remains sticking, i.e. $x(t) \in S_0$, as long as
\begin{equation}
x_1(t) = x_1(t^{s}) + x_2(t^{s}) \int \limits_{t^s}^{t^c} dt
\label{eq:sticking}
\end{equation}
does not violate \eqref{eq:theorem2}. That means, the state
trajectories will always leave $S_0$ at some $t^s < t^c < \infty$,
and always on one of the bounding straight-lines (cf. Figure
\ref{fig:2}) which are given by $| b x_2 + c x_1 | = \gamma$.
\end{rmk}

\begin{rmk}
\label{rem:5} The integral control parameter has the following
impact: (i) if $c \rightarrow 0^{+}$, the bounding straight-lines
tend to align parallel to the $x_1$-axis and, thereupon,
$(t^c-t^s) \rightarrow \infty$. This is in accord with the
analysis performed for a second-order control system with relay
\eqref{eq:relay} and without integral feedback, cf.
\cite{alvarez2000,ruderman2023}. (ii) if $c \rightarrow \infty$,
the bounding straight-lines tend to collapse with the $x_2$-axis
and, thereupon, $(t^c-t^s) \rightarrow 0$. Although this
eliminates the stiction region as such, it leads also to a
violation of the GAS condition \eqref{eq:GAScond}.
\end{rmk}

Worth noting is that the above description of the stiction region
$S_0$ and, in particular, Remark \ref{rem:5} are well in accord
with the sliding-mode dynamics once the system is on the switching
surface inside of $S_0$. Staying in the sliding-mode
(correspondingly on the switching surface $S = 0$) requires
\begin{equation}\label{eq:stayingInSM}
    \dot{S} = CA x+ CB u = 0 \quad \hbox{for all} \quad t^{s} < t <
    t^{c}.
\end{equation}
Solving \eqref{eq:stayingInSM} with respect to $u$, one obtains
the so-called equivalent control, see e.g.
\cite{shtessel2014,utkin2020}, as
\begin{equation}\label{eq:eqcontrol}
    u_{e}=-(CB)^{-1} CAx.
\end{equation}
Here we recall that an equivalent control is a linear one required
to maintain the system in an ideal sliding-mode without fast
switching. Thus, substituting \eqref{eq:eqcontrol} instead of $u$,
cf. \eqref{eq:tf}, results in an equivalent system dynamics
\begin{equation}\label{eq:eqdynamics}
    \dot{\bar{x}}=\bigl[I-B(CB)^{-1}C \bigr] A \bar{x} = \Omega \, A
    \bar{x}.
\end{equation}
Then, an equivalent state trajectory $\bar{x} = (x_1, x_2,
0)^{\top}$ is driven by \eqref{eq:eqdynamics} as long as the
system remains in sliding-mode, i.e. as long as
\eqref{eq:theorem2} is satisfied. Also recall that $\Omega$ serves
as a projection operator of the original system dynamics, while
satisfying $C \Omega = 0$ and $ \Omega B = 0$. Evaluating
\eqref{eq:eqdynamics} reveals the equivalent system dynamics as
\begin{equation}\label{eq:eqdynamicsMatrix}
    \left(%
\begin{array}{c}
  \dot{x}_{1} \\
  \dot{x}_{2} \\
  \dot{x}_{3}\\
\end{array}%
\right) =
\left(%
\begin{array}{ccc}
  0 & 1 & 0 \\
  0 & 0 & 1 \\
  0 & 0 & 0 \\
\end{array}%
\right)
\left(%
\begin{array}{c}
  x_{1}(t^{s}) \\
  x_{2}(t^{s}) \\
  0 \\
\end{array}%
\right)
\end{equation}
when $x(t) \in S_0, \; \forall \;\: t^s < t < t^c$.

For any $x_2(t_i^s)$ of entering into the stiction region $S_0$,
the integral state at leaving $S_0$ is then given by
\begin{equation}\label{eq:intexisS0}
    x_1(t_i^c) = \frac{\gamma \bigl|x_2(t_i^s)\bigr|^{-1} - b}{c} \,
    x_2(t_i^s),
\end{equation}
while $x_2(t_i^c) = x_2(t_i^s)$, cf. Figure \ref{fig:2}.
Considering the Euclidian norm of the state vector $\Phi = \|x \|$
for an unperturbed system \eqref{eq:tf}, \eqref{eq:ssmatr} during
stiction, one can recognize that the norm can either increase or
decrease between two time instants $t_i^s < t_i^c$, this is
depending on the state value $\bigl(x_1(t_i^s), x_2(t_i^s), 0
\bigr)$ of entering $S_0$. The case difference is exemplified by
two stiction trajectories shown by the blue and red solid lines
and the corresponding sphere projection (circles $\Phi'$ and
$\Phi''$) in Figure \ref{fig:2}. Note that after leaving $S_0$,
the state norm is always decreasing until the next sticking,
because the system is GAS, while the corresponding shrinking
sphere is given by
\begin{equation}\label{eq:sphere}
\Phi^c = \biggl( \Bigl(x_1 \pm \frac{\gamma}{c}  \Bigr)^2 + x_2^2
+ x_3^2 \biggr)^{\frac{1}{2}}.
\end{equation}
Note that the coordinate $(\pm \gamma/c, 0, 0)$ of the center of
the sphere \eqref{eq:sphere} results from both attraction points
of the equilibria, i.e. $\dot{x}=0 \: \Rightarrow \: c x_1 = \pm
\gamma$, cf. \eqref{eq:ode}, \eqref{eq:tf}, \eqref{eq:ssmatr}.

Against the above background, we are now to analyze the possible
scenarios of the asymptotic convergence of $x(t)$, including the
appearance of so-called \emph{stick-slip cycles}.

\subsection{Analysis of stick-slip convergence}
\label{sec:2:sub:3}

While the finite set of equilibria
\begin{equation}\label{eq:attractset}
\bar{S} \equiv \biggl\{ x \in \mathbb{R}^3 \: : \; x_1 \in
\Bigl\{-\frac{\gamma}{c}, \frac{\gamma}{c} \Bigr\},\, x_2 = 0 ,\,
x_3 = 0 \biggr\}.
\end{equation}
is globally asymptotically attractive (cf. with
\cite{bisoffi2017}), the convergence to $\bar{S}$ can appear with
an infinite number of the stick-slip cycles or, otherwise,
exponentially after experiencing at least one sticking phase, cf.
\cite{armstrong1996,bisoffi2017,ruderman2022}. The former
convergence mode is characterized by an alternating entering and
leaving $S_0$. Closer to $\bar{S}$ each next stiction trajectory
$x(t)$ proceeds, longer period $t_i^c - t_i^s$ the corresponding
sticking phase has. That means $|x_2(t_i^s)| \searrow \quad
\Rightarrow \quad (t_i^c - t_i^s) \nearrow$, cf.
\eqref{eq:sticking}, \eqref{eq:intexisS0}. While the analytic
solutions of \eqref{eq:ode} are available for both sub-spaces
$\mathbb{R}^{3}_{+}$ and $\mathbb{R}^{3}_{-}$, which are divided
by $S$, cf. \eqref{eq:slidingsurf}, see Appendix \emph{B}, their
exact analytic form is largely depending on the set of $\{a, b, c,
\gamma, C_1, C_2, C_3\}$. Especially, the initial conditions
$(C_1, C_2, C_3)$ for $t=0$ or $t=t_i^s$ can evoke the appearance
of the stick-slip cycles, examples of which will be shown in
different numerical scenarios provided in section \ref{sec:3}.

\begin{figure}[!h]
\centering
\includegraphics[width=0.9\columnwidth]{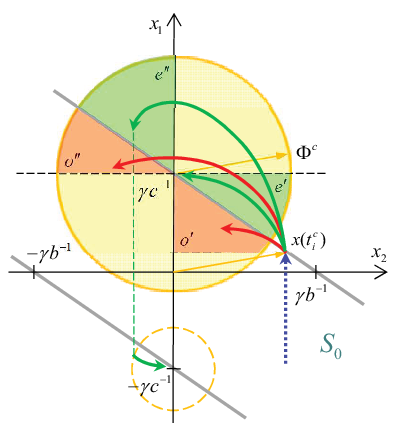}
\caption{Slipping phases after leaving $S_0$.} \label{fig:3}
\end{figure}

In the following, let us examine (qualitatively) a particular case
of the state trajectories during the slipping phase after leaving
$S_0$. That means the system was at least once in the sticking
phase, and the initial conditions $(C_1, C_2, C_3) =
\bigl(x_1(t_i^c), x_2(t_i^c), 0\bigr)$ are unambiguously given,
see \eqref{eq:intexisS0}. We consider for instance $x \in
\mathbb{R}^{3}_{-}$, this without loss of generality since the
trajectories are symmetrical with respect to the origin for $x_3 <
0$ and $x_3 > 0$. Then, the slipping trajectories are limited by a
ball $\mathcal{B} < \Phi^c$, cf. \eqref{eq:sphere}, as shown by
the $(x_1,x_2)$-projection of the ball in Figure \ref{fig:3}. Note
that for the times $t_i^c \leq t < \{t_{i+1}^s \, \vee \, \infty
\}$\footnote{The upper-limit of the time of a slipping trajectory
depends on whether a zero-crossing appears in $x(t) \in S_0$ at
some $t_i^c < t < \infty$.}, the attractor of all possible
trajectories is $\bar{S}_{+} = (\gamma c^{-1}, 0, 0)$, cf.
\eqref{eq:attractset} and Figure \ref{fig:3}.
The radius of the sphere \eqref{eq:sphere} which is enclosing the
ball
$$
\| (x_1-\gamma c^{-1}, x_2, x_3) \| < \Phi^c
$$
is determined by $\Phi^c = \| x(t_i^c) \|$, and is located in the
point where the trajectory leaves $S_0$, see the orange vector
drawn in Figure \ref{fig:3}. It is evident that for any trajectory
outgoing from $x(t_i^c)$, the $x_1(t)$-value is growing as long as
$x_2 > 0$, and then it is falling once $x_2 < 0$. Just as obvious
is that $x_2(t)$-value is decreasing only if $x(t) \in
\mathbb{R}^{3}_{-}$, the case we are considering, and it can start
to increase only after the trajectory is zero-crossing of $x_3$.
Furthermore, we note that even though $x(t_i^c)$ lies outside the
ball $\Phi^c$, a trajectory $x(t)$ is always attracted towards the
ball $\mathcal{B}$. This is due to the fact that after leaving
$S_0$ the trajectory has $\dot{x}_3 < 0$ while being driven by the
energy provided by the $(x_1,x_2)$ vector and $\Gamma < 0$, cf.
\eqref{eq:ode1} in Appendix. Thus, no trajectory can escape from
the ball $\mathcal{B}$ when the system \eqref{eq:ode1} (given in
Appendix) is GAS.

Whether the system comes again to a stiction phase, thus producing
a stick-slip cycle, depends on whether the trajectory has
zero-crossing of $x_3$ inside of the one of so-called
`oscillation' regions denoted by $o'$ or $o''$ and red-colored in
Figure \ref{fig:3}. Then, $x(t_{i+1}^s) \in S_0$ and a new
stick-slip cycle sets on. Note that whether a trajectory
experiences an overshoot and comes into the $o''$-region depends
also on the $\gamma$ and initial state values. Hence, also a
non-oscillatory linear sub-dynamics, i.e. when all roots of the
characteristic equation \eqref{eq:charpol} are real, can lead to
the trajectories landing into the $o$-regions, see numerical
examples provided in section \ref{sec:3}. On the contrary, if a
trajectory $x(t)$ does not leave one of the so-called
`exponential' regions, denoted by $e'$ or $e''$ and green-colored
in Figure \ref{fig:3}, then $x(t) \in \mathbb{R}^3 \backslash S_0$
and there is no stick-slip cycle. Also we note that if the
trajectory is overshooting into $e''$, then there is a
zero-crossing of $x_3$ and the state attractor changes unavoidably
to $\bar{S}_{-} = (-\gamma c^{-1}, 0, 0)$, cf. Figure \ref{fig:3}.

Following numerical examples disclose various asymptotic
convergence cases with and without stick-slip cycles.

%%%%%%%%%%%%%%%%%%%%%%%%%%%%%%%%%%%%%%%%%%%%%%%%%%%%%%%%%%%%%%%%%%%%%%%%%%%%%%%%
\section{Numerical examples}
\label{sec:3}

The numerical simulations of the system \eqref{eq:ode},
\eqref{eq:relay} are performed by using the state-space
realization \eqref{eq:tf}, \eqref{eq:ssmatr} and a standard
discrete-time numerical solver of MATLAB. The sampling time is set
to 0.001 sec. Note that during the slipping phases (i.e. for all
times between $t_i^c$ and $t_{i+1}^s$) the explicit solutions of
trajectories, see Appendix \emph{B}, are used. During the sticking
phases (i.e. for all times between $t_i^s$ and $t_{i}^c$), the
trajectories are computed by using the equivalent system dynamics
\eqref{eq:eqdynamicsMatrix}. This way, all trajectories are
consistent and there is no solver-related issues for simulating
the sign discontinuity \eqref{eq:relay}. Three principally
different configurations of the system parameters $a,b,c,\gamma$
are used. The first two parameter sets are the same as used in
\cite[section~IV]{bisoffi2017}, one with three distinguished
roots, and another one with one real (faster) root and two
conjugate complex roots close to the origin. The third parameter
set is assigned to be similar to the second one (i.e. with one
real and two conjugate complex roots), but with principal
difference that the real root becomes dominant, i.e. is closer to
zero. The fourth parameter set, which is similar to Example 3 in
\cite{ruderman2022}, is assigned so that to allow the slipping
phase to end in the $o'$ region, cf. Figure \ref{fig:3}.

\subsection{Three distinct real roots as in \cite{bisoffi2017}} \label{sec:3:sub:1}

\begin{figure}[!h]
\centering
\includegraphics[width=0.98\columnwidth]{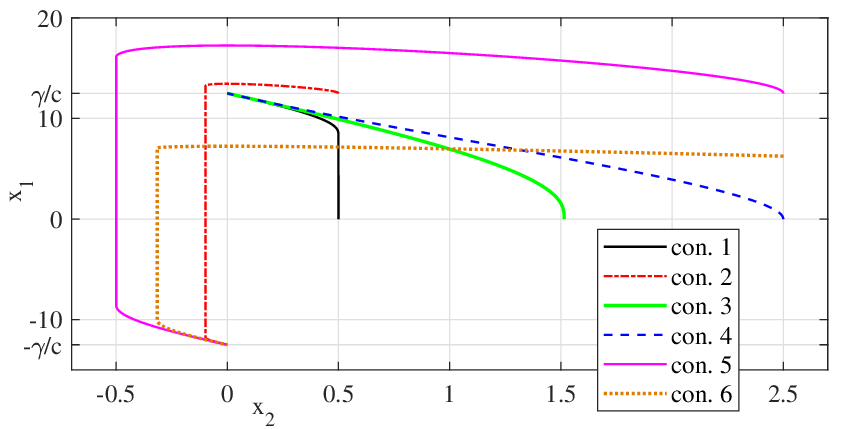}
\caption{Phase portraits in the $(x_2,x_1)$ coordinates.}
\label{fig:4}
\end{figure}

\begin{figure}[!h]
\centering
\includegraphics[width=0.98\columnwidth]{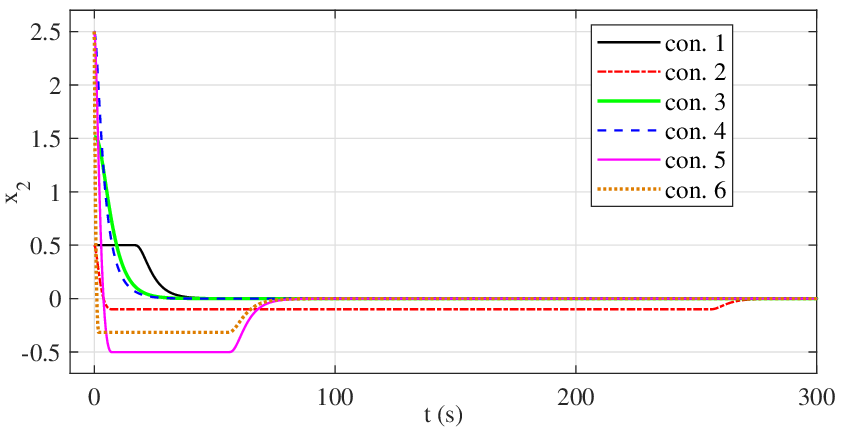}
\caption{Time series of the $x_2$ state.} \label{fig:5}
\end{figure}

Three distinct real roots $\lambda_{1,2,3} = \{-0.2, -0.5, -0.8
\}$ are assigned that correspond to $a=1.5$, $b=0.66$, $c=0.08$,
and the relay gain is $\gamma = 1$, cf. \cite{bisoffi2017}.
Different initial conditions given in Table \ref{tab:1} are used.
\begin{table}[!h]
  \renewcommand{\arraystretch}{1.4}
  \caption{Initial conditions}
  \label{tab:1}
  \footnotesize
  \begin{center}
  \begin{tabular} {|p{0.5cm}||p{0.8cm}|p{0.8cm}|p{0.8cm}|p{0.8cm}|p{0.8cm}|p{0.8cm}|}
  \hline
                     &  con. 1 & con. 2 & con. 3 & con. 4 & con. 5 & con. 6\\
  \hline \hline
  $C_1$              & 0    &  12.5    & 0       & 0      & 12.5   & 6.25      \\
  \hline
  $C_2$              & 0.5  &  0.5     & 1.5152  & 2.5    & 2.5    & 2.5     \\
  \hline
  $C_3$              & 0    &  0       & 0       & 0      & 0      &  -4.5    \\
  \hline
  \end{tabular}
  \end{center}
\end{table}
Note that the boundary $\gamma / b = 1.5152$ has also been use,
see con. 3. Also, for one set of the initial conditions (con. 1)
$x(0) \in S_0$. The phase portraits in the $(x_2,x_1)$ coordinates
are shown in Figure \ref{fig:4}. One can notice that for different
initial conditions, the trajectories experience one stick-slip
cycle, after which they converge exponentially to one of the
attraction points $(\pm \gamma c^{-1},0,0)$.
The corresponding time series of the $x_2$ state of interest
(relative displacement in case of the motion control) are shown in
Figure \ref{fig:5}. One can recognize the largely varying period
of the sticking phase, depending on $x(t^s)$ where the trajectory
entered $S_0$.

\subsection{One real and two conjugate complex roots} \label{sec:3:sub:2}

One real root and two conjugate complex roots $\lambda_{1,2,3} =
\{-6.01, -0.19 \pm j 0.79 \}$ are assigned, that corresponds to
$a=6.4$, $b=3$, $c=4$, and the relay gain is $\gamma = 1$, cf.
\cite{bisoffi2017}. The pair of the conjugate complex roots
corresponds to the eigenfrequency $\omega_0 = 0.8125$ and damping
ratio $\delta = 0.2338$, cf. Appendix \emph{B}. One can recognize
that the conjugate complex roots clearly dominate over the real
one, since lying much closer to the origin in the complex plane.

In addition, another configuration of the conjugate complex roots
with $\omega_0 = 2$ and the same $\delta = 0.2338$ was assigned.
Here the real root was selected to be clearly dominating over the
conjugate complex pair, thus resulting in $\lambda_{1,2,3} =
\{-0.2, -0.468 \pm j 1.94 \}$. The relay gain has the same
$\gamma=1$ value, while the resulting coefficients are $a=1.135$,
$b=4.187$, $c=0.8$.

\begin{figure}[!h]
\centering
\includegraphics[width=0.98\columnwidth]{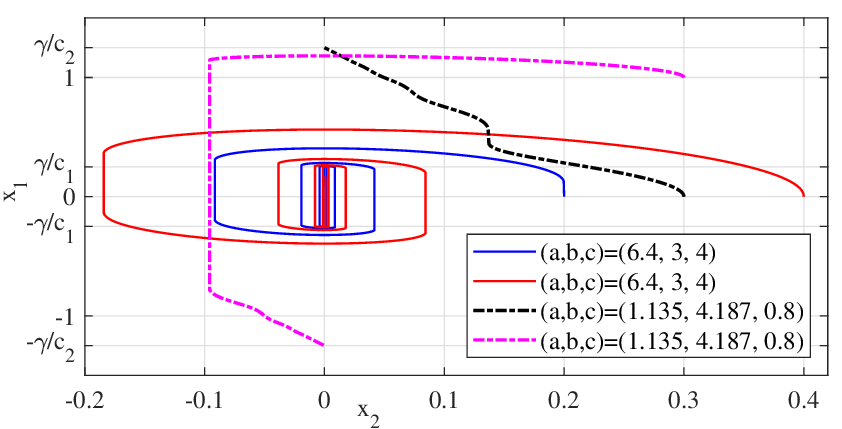}
\caption{Phase portraits in the $(x_2,x_1)$ coordinates.}
\label{fig:6}
\end{figure}

\begin{figure}[!h]
\centering
\includegraphics[width=0.98\columnwidth]{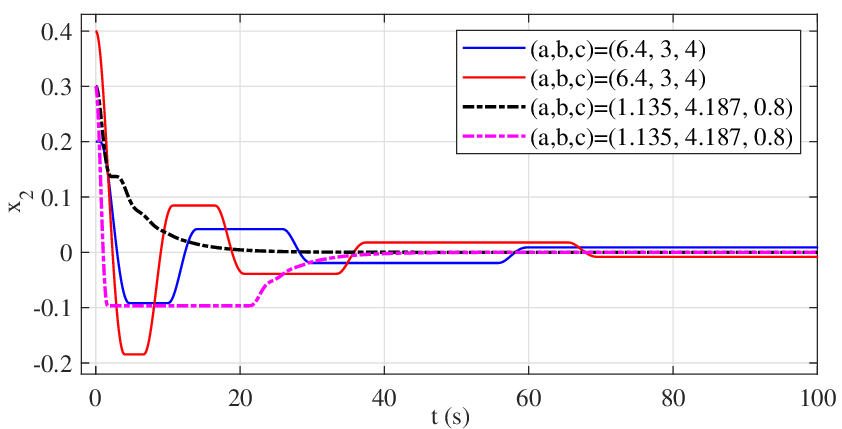}
\caption{Time series of the $x_2$ state.} \label{fig:7}
\end{figure}

For both configurations of the parameters, the phase portraits are
shown in Figure \ref{fig:6} in the $(x_2,x_1)$ coordinates, and
the corresponding time series of the $x_2$ state are shown in
Figure \ref{fig:7}.
The first configuration of parameters is depicted by the solid
lines, with two sets of the initial values $(C_1,C_2,C_3) =
\{(0,0.2,0), \: (0,0.4,0) \}$, one starting inside and other
outside of $S_0$. The second configuration of parameters is
depicted by the dash-dotted lines, with two sets of the initial
values $(C_1,C_2,C_3) = \{(0,0.3,0), \: (1,0.3,0) \}$, both
starting outside of $S_0$.
One can recognize that the first configuration of parameters (i.e.
with the dominant conjugate complex roots) comes always to the
sustained stick-slip cycles, thus circulating around both
attraction points $(\pm \gamma c_1^{-1},0,0)$. On the contrary,
the second configuration of parameters (i.e. with the dominant
real root), although also oscillatory, comes to only one
stick-slip cycle, after which the trajectories are converging
exponentially to one of the attraction points $(\pm \gamma
c_2^{-1},0,0)$.

\subsection{Conjugate complex roots as in \cite{ruderman2022}} \label{sec:3:sub:3}

One real root and two conjugate complex roots $\lambda_{1,2,3} =
\{-7.82, -1.09 \pm j 32 \}$ are assigned that correspond to
$a=10$, $b=1040$, $c=8000$, and the relay gain is $\gamma = 100$,
cf. \cite{ruderman2022}. The pair of the conjugate complex roots
corresponds to the eigenfrequency $\omega_0 = 32$ and damping
ratio $\delta = 0.034$. One can recognize that while the conjugate
complex roots dominate over the real one, the otherwise highly
oscillatory (with respect to $\delta$) system is mostly damped by
the nonlinear relay action. Thus, it is expected to come
repeatedly into $S_0$ within $o'$ region, i.e. in a close vicinity
to $x_2 = 0$, cf. Figure \ref{fig:3}.
\begin{figure}[!h]
\centering
\includegraphics[width=0.98\columnwidth]{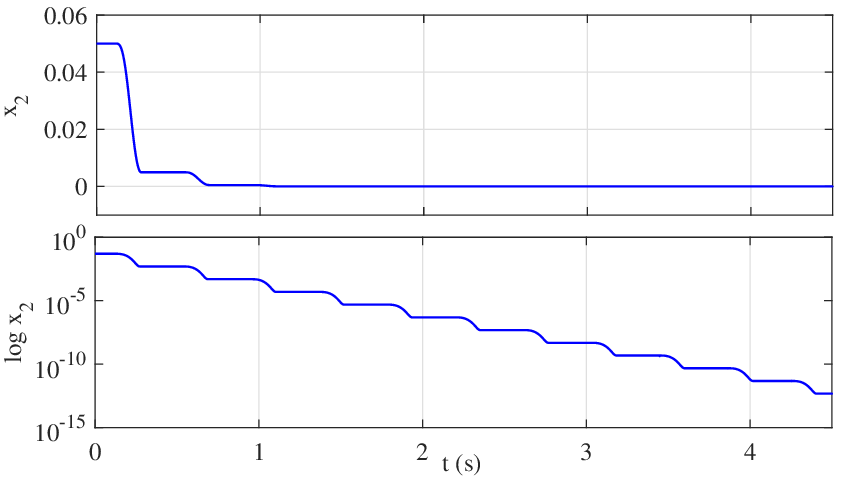}
\caption{Time series of the $x_2$ state.} \label{fig:8}
\end{figure}
The time series of the $x_2$ state are shown in Figure
\ref{fig:8}, once on the linear and once on the logarithmic scale.

\section{Experimental evidence}
\label{sec:4}

To (qualitatively) support the
results and conclusions of the analysis conducted, a series of
laboratory experiments of the motion control are exemplary shown.
The linear displacement actuator with one degree of freedom is
feedback regulated by a standard proportional-integral (PI)
control
\begin{equation}\label{eq:PIcontrol}
u(t) = K_p \, e(t) + K_i \int e(t) \, dt.
\end{equation}
Here, the (unit-less) input control force is $u(t)$, and the
output feedback control error is $e(t) = r(t) - x_2(t)$. The
reference value is assigned to be a set positioning constant $r(t)
= 0.01$ m, while the output displacement of interest $x_2(t)$ is
provided by a contactless sensor, which exhibits a relatively high
level of the measurement noise. For a detailed description of the
used actuator system, an interested reader is referred to
\cite{ruderman2022motion}.

Worth noting is that the control
gains $K_p, K_i > 0$ directly enter into the coefficients $b$ and
$c$, correspondingly, of the third-order closed-loop dynamic
system, cf. \eqref{eq:ode}. Also understood is that a weakly known
(and even uncertain) Coulomb friction coefficient of the motion
system, normalized by its lumped mass parameter, constitutes the
relay gain factor $\gamma$, cf. \eqref{eq:ode}.

\begin{figure}[!h] \centering
\includegraphics[width=0.98\columnwidth]{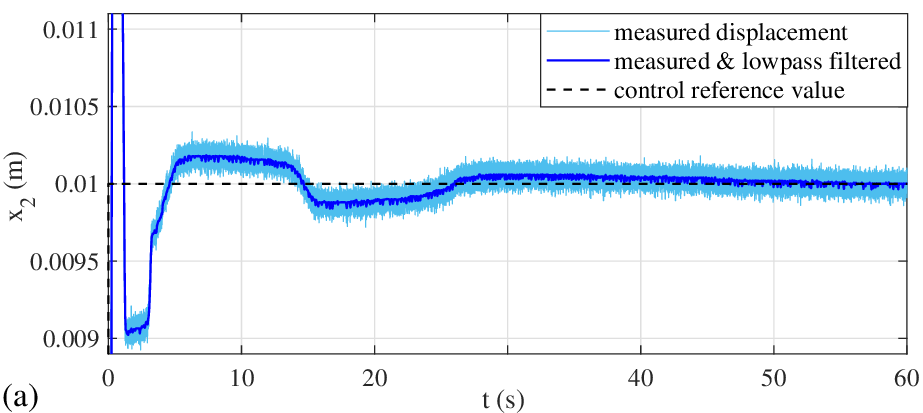}
\includegraphics[width=0.98\columnwidth]{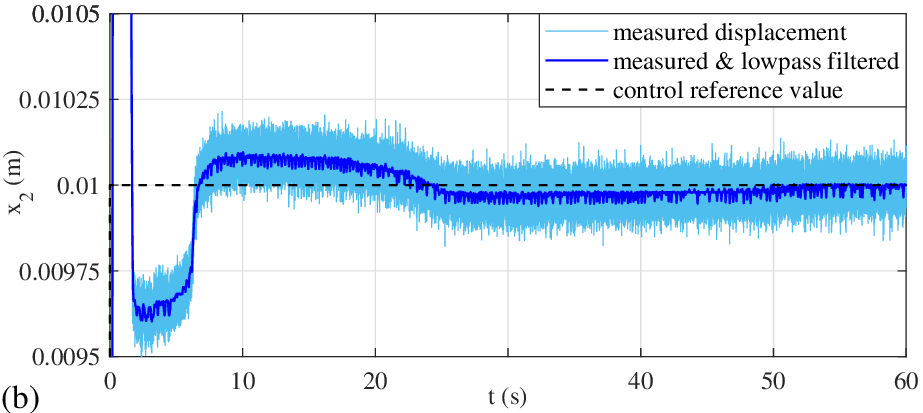}
\caption{Time series of the
experimentally measured $x_2$ state of the PI feedback controlled
linear actuator displacement: (a) control gains $K_p=50$,
$K_i=400$, (b) control gains $K_p=100$, $K_i=400$.} \label{fig:9}
\end{figure}
Two sets of the control parameters
are evaluated, as shown in Figure \ref{fig:9}, $K_p=50$, $K_i=400$
depicted in (a), and $K_p=100$, $K_i=400$ depicted in (b). For
better clarity, in addition to the raw sensor measurement, the
same output displacement $x_2(t)$, but low-pass filtered (with 10
Hz cutoff frequency), is also shown. One can recognize that in
both cases the $x_2(t)$ trajectory converges with stick-slip
cycles, while their number and period differ. It is also visible
that the sticking phase does not have an ideal horizontal plateau,
cf. with Figs. \ref{fig:5} and \ref{fig:7}. This is natural since
the real Coulomb friction is combined (correspondingly coupled)
with viscous, adhesive, creeping and other nonlinear by-effects,
cf. e.g. \cite{armstrong1994}, and the frictional contact surfaces
in the actuator system are highly inhomogeneous. Nevertheless, a
principal pattern of the alternating sticking and slipping phases,
with an increasing period for the decreasing $|e(t)|$, is in line
with the developed analysis and (ideal) numerical simulations.

\section{Summary and discussion}
\label{sec:5}

In this work, a relatively simple (to follow) analysis of the
convergence behavior of the third-order dynamic systems which have
integral feedback action and discontinuous relay disturbance was
provided. The GAS property of the system was proven and the
so-called stiction regime in a well-defined region of the
state-space was specified and analyzed by using the sliding-mode
principles. The corresponding equivalent dynamics was established,
and the conditions for entering into and leaving from the stiction
region were described. The possible scenarios of a stick-slip
convergence were addressed and illustrated by the purposefully
designed dedicated numerical examples. Moreover, an illustrative
experimental example of a PID feedback controlled positioning in
presence of the Coulomb friction of an actuator was demonstrated
in favor of the developed and discussed analysis of the
convergence.

The following can be summarized and brought into discussion. For
the real roots of characteristic polynomial of the linear
subsystem, only one stick-slip cycle can appear (depending on the
initial conditions) before the exponential convergence takes place
towards one of the two symmetrical points of attraction. For the
conjugate complex roots, the appearance of persistent stick-slip
cycles can depend on whether the conjugate complex roots dominate
over the real one. If such system is predominantly damped by the
relay feedback action, the persistent stick-slip cycles can appear
even without overshoot, see example shown in section
\ref{sec:3:sub:3}. If, otherwise, the real root dominates over the
conjugate complex pair, certain low number of the stick-slip
cycles (but not a persistent series) is expected before entering
into the exponential convergence. However, such statements made
above require a more detailed analysis of both, the polynomial
coefficients and the initial conditions, while taking into account
also the relay gain value. The latter enters directly into the
particular solutions shown in Appendix \emph{B}. At large, the
present work, in addition to the previously published works
\cite{armstrong1996,bisoffi2017,ruderman2022}, provide further
insight into the problem of state feedback plus integral control
behavior of the motion systems in presence of the discontinuous
Coulomb friction. It might also be interesting to consider order higher than three in the future for
a class of systems with feedback relay disturbance in extension of
\eqref{eq:ode}.

\section*{Appendix}
\label{sec:5}

\subsection*{A: Solution of Lyapunov equation} \label{sec:5:sub:1}

For $Q=I$ and $A$ given in \eqref{eq:ssmatr}, the unique solution
of \eqref{eq:lyapeq} is the $3 \times 3$ symmetric matrix $P$ with
the elements
\begin{eqnarray}
\nonumber  p_{1,1} &=& -\frac{c(a^2+ac-b+c^2)+ab^2}{2c(c-ab)}, \\
\nonumber  p_{1,2} &=& -\frac{a^2b+c^2(b+1)}{2c(c-ab)}, \\
\nonumber  p_{1,3} &=& \frac{1}{2c}, \\
\nonumber  p_{2,2} &=& -\frac{a(a^2+ac+c^2)+c(b^2+b+1)}{2c(c-ab)}, \\
\nonumber  p_{2,3} &=& -\frac{a(a+c)+c^2}{2c(c-ab)}, \\
\nonumber  p_{3,3} &=& -\frac{a+c(b+1)}{2c(c-ab)}.
\end{eqnarray}

\subsection*{B: Particular solutions during the system slipping} \label{sec:5:sub:2}

For the system \eqref{eq:ode}, \eqref{eq:relay}, equivalently
\eqref{eq:tf}, \eqref{eq:ssmatr}, the particular solutions during
slipping, i.e. for $\ddot{y}=x_3\neq0$, are those of the following
nonhomogeneous differential equation
\begin{equation}
\begin{split}
\dot{x}_3(t) + a \, x_3(t) + b \, x_2(t) + c \, x_1(t) = - \Gamma
,
\\[0.5mm]
\qquad \forall \; \; t_i^c \leq t  \, \left\{%
\begin{array}{ll}
    \leq t_0       & \hbox{if } x_3(t_0) = 0, \\
    < \infty       & \hbox{if } x_3(t) \neq 0,\\
\end{array}%
\right. \end{split} \label{eq:ode1}
\end{equation}
where the constant exogenous right-hand-side $\Gamma \equiv \gamma
\, \mathrm{sign}(x_3)$ of \eqref{eq:ode1} is assigned for the
particular trajectories $x \in \mathbb{R}^{3}_{-} \vee
\mathbb{R}^{3}_{+}$. Next, we need to distinguish whether the
corresponding characteristic polynomial \eqref{eq:charpol} has all
real roots or a pair of the conjugate complex roots.

\subsubsection*{Three distinct real roots}

Assume that \eqref{eq:charpol} has three distinct real roots
$\lambda_1 \neq \lambda_2 \neq \lambda_3$ while the dynamic system
\eqref{eq:ode1} is asymptotically stable, i.e. $\lambda_1,
\lambda_2, \lambda_3 > 0$ and the inequality \eqref{eq:GAScond}
holds for
\begin{equation}
a = \lambda_1 + \lambda_2 + \lambda_3, \; b = \lambda_2 \lambda_3
+ \lambda_1 \lambda_2 + \lambda_1 \lambda_3, \; c = \lambda_1
\lambda_2 \lambda_3. \label{eq:coeffic1}
\end{equation}
Then, the particular solution of \eqref{eq:ode1} is given by
\begin{equation}
x_1(t) = -\frac{\Gamma}{\lambda_1 \lambda_2 \lambda_3} + K_1
e^{-\lambda_1 t} - K_2 e^{-\lambda_2 t} + K_3 e^{-\lambda_3 t},
\label{eq:ode1sol1}
\end{equation}
with the coefficients
\begin{eqnarray}
\nonumber K_1 &=& \frac{\Gamma + C_3 \lambda_1 + C_2 \lambda_1
(\lambda_2 + \lambda_3) + C_1 \lambda_1 \lambda_2 \lambda_3
}{\lambda_1(\lambda_1-\lambda_2)(\lambda_1-\lambda_3)},
\\[0.5mm]
\nonumber K_2 &=& \frac{\Gamma + C_3 \lambda_2 + C_2 \lambda_2
(\lambda_1 + \lambda_3) + C_1 \lambda_1 \lambda_2
\lambda_3}{\lambda_2(\lambda_1-\lambda_2)(\lambda_2-\lambda_3)},
\\[0.5mm]
\nonumber K_3 &=& \frac{\Gamma + C_3 \lambda_3 + C_2 \lambda_3
(\lambda_1 + \lambda_2) + C_1 \lambda_1 \lambda_2
\lambda_3}{\lambda_3(\lambda_1-\lambda_3)(\lambda_2-\lambda_3)}.
\label{eq:coeffOde1}
\end{eqnarray}
Respectively,
\begin{equation}
x_2(t) = \dot{x}_1(t), \quad  x_3(t) = \ddot{x}_1(t).
\label{eq:derivOde1}
\end{equation}

\subsubsection*{One real and two conjugate complex roots}

Assume that \eqref{eq:charpol} has one real root $\lambda_1$ and a
pair of conjugate complex roots $\lambda_{2,3} = -\delta \omega_0
\pm j \omega_0 \sqrt{1-\delta^2}$ which are parameterized by the
eigenfrequency $\omega_0 > 0$ and damping ratio $0 < \delta \leq
1$. Recall that the case of the double real roots $\lambda_2 =
\lambda_3$ is covered by $\delta = 1$. The dynamic system
\eqref{eq:ode1} is asymptotically stable, i.e. $\lambda_1,
\lambda_2, \lambda_3 > 0$ and the inequality \eqref{eq:GAScond}
holds for
\begin{equation}
a = 2\delta \omega_0 + \lambda_1, \;\; b = 2\delta \omega_0
\lambda_1 + \omega_0^2, \;\; c = \lambda_1 \omega_0^2.
\label{eq:coeffic2}
\end{equation}
Then, the particular solution of \eqref{eq:ode1} is given by
\begin{eqnarray}
\label{eq:ode1sol2}
x_1(t) & = & -\frac{\gamma }{\lambda _{1}\,{\omega _{0}}^2} \\[0.5mm]
\nonumber & + & \frac{{\mathrm{e}}^{-\lambda
_{1}\,t}\,\left(C_{1}\,\lambda _{1}\,{\omega
_{0}}^2+2\,C_{2}\,\delta \,\lambda _{1}\,\omega _{0}+\gamma
+C_{3}\,\lambda _{1}\right)}{\lambda _{1}\,\left({\lambda
_{1}}^2-2\,\delta \,\lambda _{1}\,\omega _{0}+{\omega
_{0}}^2\right)} \\[0.5mm]
\nonumber & + & \frac{K_4 {\mathrm{e}}^{-\omega
_{0}\,t\,\left(\delta -\sqrt{\delta ^2-1}\right)}}{2\,{\omega
_{0}}^2\,\sqrt{1-\delta }\,\sqrt{-\delta -1}\,\left(\lambda
_{1}-\delta \,\omega
_{0}+\omega _{0}\,\sqrt{\delta ^2-1}\right)} \\[0.5mm]
\nonumber & + & \frac{K_5{\mathrm{e}}^{-\omega
_{0}\,t\,\left(\delta +\sqrt{1-\delta }\,\sqrt{-\delta
-1}\right)}}{2\,{\omega _{0}}^2\,\sqrt{1-\delta }\,\sqrt{-\delta
-1}\,\left(\delta \,\omega _{0}-\lambda _{1}+\omega
_{0}\,\sqrt{\delta ^2-1}\right)},
\end{eqnarray}
with the coefficients
\begin{eqnarray}
\nonumber K_4 &=& C_{3}\,\omega _{0}+\delta \,\gamma +\gamma
\,\sqrt{\delta ^2-1}+C_{2}\,\lambda _{1}\,\omega
_{0}+C_{2}\,\delta \,{\omega _{0}}^2 \\[0.5mm]
\nonumber & + & C_{2}\,{\omega _{0}}^2\,\sqrt{\delta ^2-1} + C_{1}
\, \lambda_{1} \, {\omega _{0}}^2 \Bigl(\delta + \sqrt{\delta
^2-1} \Bigr),
\\[0.5mm]
\nonumber K_5 &=& C_{3}\,\omega _{0}+\delta \,\gamma -\gamma
\,\sqrt{\delta ^2-1}+C_{2}\,\lambda _{1}\,\omega
_{0}+C_{2}\,\delta \,{\omega _{0}}^2 \\[0.5mm]
\nonumber & - & C_{2}\,{\omega _{0}}^2\,\sqrt{\delta ^2-1} + C_{1}
\, \lambda_{1} \, {\omega _{0}}^2 \Bigl(\delta - \sqrt{\delta
^2-1} \Bigr). \label{eq:coeffOde2}
\end{eqnarray}
Respectively, \eqref{eq:derivOde1} applies as well.

\bibliographystyle{elsarticle-num}
\bibliography{references}             % bib file to produce the bibliography
                                                     % with bibtex (preferred)

\end{document}